\documentclass[12pt]{article}
\usepackage{amsmath}
\usepackage{amssymb, amsbsy, amstext}
\usepackage{srcltx}
\usepackage{amscd,amsfonts,color}
\usepackage{latexsym}
\usepackage{epic,eepic,epsfig}
\usepackage{tikz}
\usepackage{enumerate}

\usetikzlibrary{decorations.markings}
\tikzstyle{vertex}=[circle, draw, inner sep=0pt, minimum size=5pt]
\tikzstyle{vedge}=[circle, draw, inner sep=0pt, minimum size=5pt]
\newcommand{\vertex}{\node[vertex]}
\newcommand{\vedge}{\node[vedge]}
\usepackage{pgflibraryarrows}
\usepackage{pgflibrarysnakes}

\input amssym.def
\input amssym.tex

\newtheorem{lem}{Lemma}[section]%
\newtheorem{theo}[lem]{Theorem}%
\newtheorem{alm}[lem]{Algorithm}%
\newtheorem{cor}[lem]{Corollary}%
\newtheorem{prop}[lem]{Proposition}%
\newtheorem{rem}[lem]{Remark}%
\newtheorem{defi}[lem]{Definition}%
\newtheorem{exa}[lem]{Example}%

\newcommand{\B}[1]{\mathbb #1}  \newcommand{\V}{\hbox{\rm V}}  \def\mod{\hbox{\rm mod }}
\newcommand{\demo}{{\bf Proof}\hskip10pt}   \newcommand{\qed}{\hfill $\Box$}
\newcommand{\ox}{\overline x}   
\newcommand{\oa}{\overline a}

\newcommand{\Sym}{\hbox{\rm Sym}}
  \newcommand{\Mon}{\hbox{\rm Mon}}     
  \newcommand{\Aut}{\hbox{\rm Aut}}     
   \newcommand{\Core}{\hbox{\rm Core}}  \newcommand{\Syl}{\hbox{\rm Syl}}

\newcommand{\AGL}{\hbox{\rm AGL}}  \newcommand{\PSU}{\hbox{\rm PSU}}  \newcommand{\PSp}{\hbox{\rm PSp}}
\newcommand{\AG}{\hbox{\rm AG}}    \newcommand{\PG}{\hbox{\rm PG}}    
\newcommand{\GL}{\hbox{\rm GL}}      \newcommand{\PSL}{\hbox{\rm PSL}}
\newcommand{\G}{\mathcal{G}}

\def\a{\alpha}     \def\th{\theta} \def\g{\gamma}    
  \def\di{\bigm|} \def\lg{\langle}  \def\rg{\rangle}
\def\HH{\mathcal{H}}

\begin{document}

\begin{center} {\bf\large Regular embedding of simple hypergraphs}
\footnote{Corresponding to pktide@163.com.
E-mail addresses:    zhuyanhong911@163.com (Yanhong Zhu), pktide@163.com (Kai Yuan).}
\bigskip

\end{center}

\begin{center}  Yanhong Zhu$^a$ and   Kai Yuan$^{b}$
\end{center}

\begin{center}
{\footnotesize  $^a$School of Mathematical Sciences, Liaocheng University, Liaocheng 252000, People's Republic of China}
\\{\footnotesize  $^{b}$School of Mathematics and Information Sciences, Yantai University, Yantai 264005,  People's Republic of China}
\end{center}

\vskip 5mm

\begin{abstract} Regular hypermaps with underlying simple hypergraphs
are analysed. We obtain an algorithm to classify the regular embeddings of simple hypergraphs with given
order, and determine the automorphism groups of
regular embedding of simple hypergraphs with prime square order.

\vskip 2mm

\medskip


\end{abstract}

\medskip

\section{Introduction}

Let $\varGamma$ be a hypergraph. The associated {\it Levi graph} $\mathcal{G}$ of $\varGamma$ is a bipartite graph with
vertex set $V(\varGamma) \sqcup E(\varGamma)$, in which $v \in V(\varGamma)$ and $e \in E(\varGamma)$
are adjacent if and only if $v$ and $e$ are incident in $\varGamma$.  The hypergraph $\varGamma$ is called  simple  if for any hyperedge $e$, $e$ is not a subset of some other hyperedges.

Topologically, a {\it hypermap} is a $2$-cell embedding  of the Levi graph $\mathcal{G}$ of a connected hypergraph $\varGamma$ (may have multiple edges) into a compact and connected surface $\mathcal{S}$ without border,
where the vertices of $\mathcal{G}$ in two biparts are respectively called the hypervertices and hyperedges of the hypermap, and the connected regions of $\mathcal{S}\backslash\mathcal{G}$ are called {\it hyperfaces}. $\varGamma$ is called the {\it underlying
hypergraph} of $\mathcal{H}$, and $\mathcal{G}$ is called the underlying Levi graph of $\mathcal{H}$.
 A hypermap $\mathcal{H}$ is {\it orientable} if the underlying surface $\mathcal{S}$ is orientable.

The automorphism group of a regular hypermap does not necessarily have a faithful action
on its  hypervertices. The aim of this paper is to give a description of the associated hypermaps with underlying simple hypergraphs and groups acting
on hypervertices. A need for such analysis is motivated by the works of regular maps with simple underlying graphs \cite{B,CM,DKN,JJ,Wil} and application of results from permutation groups to regular hypermaps.

One significant problem in the theory of hypermaps is to
characterize hypermaps with high symmetry, that is, characterize hypermaps
 of a given supporting surface \cite{C}, characterize hypermaps
 of a given underlying hypergraph \cite{JM,LK, YW2},
characterize hypermaps
 of a given underlying Levi graph \cite{CH}, characterize hypermaps
 of a given number of hyperfaces \cite{DH} and characterize hypermaps
 of a given automorphism group \cite{CPS, DY, HD, YW}. In this paper, we characterize regular hypermaps with underlying simple hypergraphs of prime square order.

This paper is  organised as follows. In Section 2, we introduce embedding of simple hypergraphs. Some preliminary group theoretical results will be given in Sections 3. The automorphism groups of regular embedding of simple hypergraphs with prime square order are determined in Section 4 and 5.

\section{Embedding of  simple hypergraphs}

\subsection{Hypergraphs}

A finite multi-subset $E$
of $X$ is a set in which finite multiplicity of elements is allowed
and taken into account.

\begin{defi}
A \textbf{hypergraph} $H$ is a pair $(V,E)$ where $V$ is a set of vertices and $E$ is a muti-subset of $\mathcal{P}(V)$ which is a set of non-empty subsets of $V$. The elements of $E$ are called hyperedges.
\end{defi}

\begin{defi}
Let $H= (V,E)$ be a hypergraph with $|V| =n$ and $|E| =m$.
Then the {\bf Levi graph} is the bipartite graph $G =
(V \sqcup E, E')$, where $(v_i,
e_j) \in E'$ if and only if $v_i \in e_j$.
\end{defi}

\begin{exa}
Let $H$ be the hypergraph with vertex set $V = \{1,2,3,4,5\}$ and edge
set
	\[ E = \{e_1=\{1,3\},e_2=\{1,3\}, e_3=\{2,3\}, e_4=\{1,3,5\}\}.	\]
 The following is the Levi graph representation of $H$.

\begin{center}
\[\begin{tikzpicture}[
	every edge/.style={
        draw,
        postaction={decorate,
                    decoration={markings,mark=at position 1 with {\arrow{>}}}
                   }
        }
        ]
	\vertex[fill] (s5) at (0,1) [label=left:${5}$] {};
	\vertex[fill] (s4) at (0,2) [label=left:${4}$] {};
	\vertex[fill] (s3) at (0,3) [label=left:${3}$] {};
	\vertex[fill] (s2) at (0,4) [label=left:${2}$] {};
	\vertex[fill] (s1) at (0,5) [label=left:${1}$] {};
	\vedge (t4) at (1,1.5) [label=right:$e_{4}$] {};
	\vedge (t3) at (1,2.5) [label=right:$e_{3}$] {};
	\vedge (t2) at (1,3.5) [label=right:$e_{2}$] {};
	\vedge (t1) at (1,4.5) [label=right:$e_{1}$] {};

	 	\draw (s1) -- (t1) (s3) -- (t1)
	 	(s1) -- (t2)
	 	(s1) -- (t4)
	 	(s2) -- (t3)
	 	(s3) -- (t2)
	 	(s3) -- (t3)
	 	(s3) -- (t4)
	 	(s5) -- (t4)
	 	
	;
\end{tikzpicture}\]

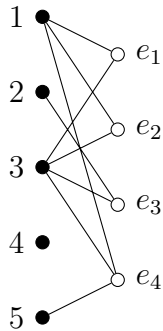
Fig. $2$:\, {\it The Levi graph representation of $H$.}
\end{center}

\end{exa}

\subsection{Hypermaps}

For a topological hypermap ${\cal H}$ with simple underlying Levi graph, choose a center for each  hyperface and  subdivide the hypermap by adjoining the hyperface centers to its adjacent hypervertices and hyperedges. Then we get a triangular subdivision whose triangles are the {\it flags} of this  hypermap,    which  are represented by little triangle around  hypervertices  as shown in Fig. 1 (right). Now, we define three involuntary permutations $\g_0$, $\g_1$ and $\g_2$  on     the set $F$  of flags of ${\cal H}$ as follows:
 $\gamma_0$  exchanges two flags adjacent to  the same hyperedge and  center but  distinct  hypervertices;
 $\gamma_1$   exchanges two flags adjacent to the same hypervertex and  center but   distinct  hyperedges;  and $\gamma_2$ exchanges two flags  adjacent  to the same hypervertex and hyperedge but  distinct centers,  see Fig. 1 (left).

\begin{center}
  \begin{tikzpicture}[
dashed line/.style={dashed, thin},
dot/.style = {
      draw,
      fill = white,
      circle,
      inner sep = 0pt,
      minimum size = 5pt
    }
]
\draw (1,2) -- (2,1) -- (3,2) (1,-2) -- (2,-1) -- (3,-2) (2,1) -- (2,-1);

\draw[dashed](0,0) -- (1,2) (0,0) -- (2,1) (0,0) -- (2,-1) (0,0) -- (1,-2) (4,0) --  (3,2) (4,0) -- (2,1) (4,0) -- (2,-1) (4,0) -- (3,-2) ;

\fill (2,1) circle (3pt); \fill(1,-2) circle (3pt); \fill(3,-2) circle (3pt);

\draw (1,2) node[dot] {} (3,2) node[dot] {} (2,-1) node[dot] {} ;

\path (1,0) node[label = {right: $f^{\gamma_2}$}] {};
\path (3,0) node[label = {left: $f$}] {};
\path (3,1.5) node[label = {below: $f^{\gamma_1}$}] {};
\path (3,-1.5) node[label = {above: $f^{\gamma_0}$}] {};

\node at (2,2) {$\cdots$};

\node at (2,-2) {$\cdots$};
  \begin{scope}[xshift=7cm]

\draw (1,2) -- (2,1) -- (3,2) (1,-2) -- (2,-1) -- (3,-2) (2,1) -- (2,-1);

\fill (2,1) circle (3pt); \fill(1,-2) circle (3pt); \fill(3,-2) circle (3pt);

\draw (1,2) node[dot] {} (3,2) node[dot] {} (2,-1) node[dot] {} ;

\fill (2.26,1.12)--(2.52,1.38)--(2.8,1.12)--cycle;
\fill (1.5,0.44)--(1.9,0.44)--(1.9,0.84)--cycle;
\fill (1.08,-1.82)--(1.08,-1.26)--(1.34,-1.52)--cycle;
\fill (2.94,-1.82)--(2.68,-1.54)--(2.96,-1.26)--cycle;
\fill (2.1,0.44)--(2.1,0.84)--(2.5,0.44)--cycle;
\fill (1.2,1.12)--(1.76,1.12)--(1.5,1.38)--cycle;

\path (1.7,0.42) node[label = {below: $f^{\gamma_2}$}] {};
\path (2.3,0.42) node[label = {below: $f$}] {};
\path (2.62,1.24) node[label = {right: $f^{\gamma_1}$}] {};
\path (2.94,-1.54) node[label = {right: $f^{\gamma_0}$}] {};
\node at (2,2) {$\cdots$};
\node at (2,-2) {$\cdots$};
  \end{scope}
\end{tikzpicture}

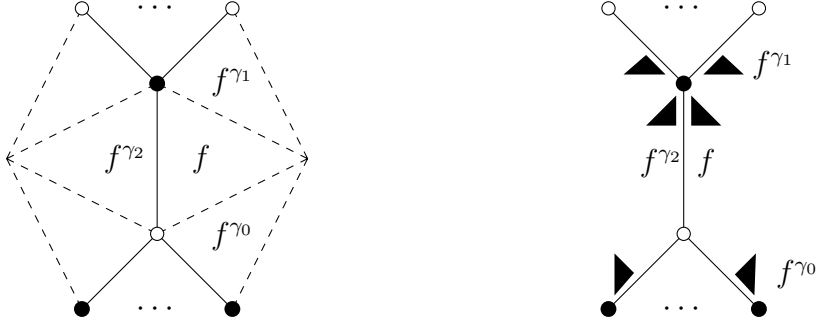
Fig. 1:\, {\it Flags are drawn in the  left figure and   represented by little triangles around  hypervertices  in the right figure.}
\end{center}

The group $\lg \g_0, \g_1, \g_2\rg $ is called the monodromy group of ${\cal H}$, denoted by $\Mon({\cal H})$. Since ${\cal G}$ is connected, $\Mon({\cal H})$ acts transitively  on $F$. We define
the hypervertices, hyperedges and hyperface boundaries of ${\cal H}$ to be the orbits of the subgroups $\lg \g_1, \g_2\rg$,$\lg \g_0, \g_2\rg$ and $\lg \g_0, \g_1\rg$, respectively. The incidence
in ${\cal H}$ can be represented by nontrivial intersection. Given two hypermaps ${\cal H}={\cal H}(F;\g_0, \g_1, \g_2)$ and ${\cal H}^\prime={\cal H}(F^\prime;\g_0^\prime, \g_1^\prime,\g_2^\prime)$, a morphism
$\sigma$ from ${\cal H}$ to ${\cal H}^\prime$ is a map $\sigma: F\rightarrow F^\prime$ such that $ \g_0\sigma=\sigma \g_0^\prime, \g_1\sigma=\sigma \g_1^\prime$ and $\g_2\sigma=\sigma \g_2^\prime$. If $\sigma$ is also bijective,
we say that $\sigma$ is an isomorphism. In particular, if ${\cal H}={\cal H}^\prime$, then $\sigma$ is called an automorphism of
${\cal H}$. The automorphisms of ${\cal H}$ form a group $\Aut({\cal H})$ which is called the automorphism group of the hypermap ${\cal H}$. By the definition
of hypermap isomorphism, we have $\Aut({\cal H})=C_{\Sym(F)}(\Mon({\cal H}))$, the centralizer of $\Mon({\cal H})$ in $\Sym(F)$ . Therefore, $\Aut({\cal H})$ acts semiregularly on $F$, noting the transitivity of $\Mon({\cal H})$ on $F$.
A hypermap is called {\it regular} if its automorphism group acts regularly on the set of all flags
(incident hypervertex-hyperedge-hyperface triples). Suppose that ${\cal H}={\cal H}(F;\g_0, \g_1, \g_2)$ is regular. This allows us to identify $\Mon({\cal H})\cong \Aut({\cal H})\cong G$ as abstract groups, though as permutation groups on $F$ they are
distinct, and can be regarded as the right and left regular representations of the same group $G$. Therefore, we may identify $F$ with $G$, so that ${\cal H} \cong {\cal H}(G;\g_0, \g_1,\g_2)$. The order $|\g_1\g_2|=k$, $|\g_2\g_0|=m$ and $|\g_0\g_1|=m$ of the elements
$\g_1\g_2$, $\g_2\g_0$ and $\g_0\g_1$ are called the {\it valency} of a hypervertex,
  hyperedge and hyperface, respectively. The triple $(k,m, n)$ is called the
{\it type} of the regular hypermap ${\cal H}$.

\vskip 0.3cm

In all group presentations in this article we will assume that the powers
at generators and relators are their {\it true orders}.

Let $G$ be an arbitrary group with partial presentation
\begin{equation}
G := \langle \g_0, \g_1, \g_2 \di \g_0^2=\g_1^2=\g_2^2=(\g_1\g_2)^k=(\g_2\g_0)^m=(\g_0\g_1)^n=\cdots = 1\rangle \label{XX}
\end{equation}
where dots indicate the (possible) presence of other relations and let $\HH = \HH(G; \g_0, \g_1, \g_2)$ be the associated regular hypermap.

\begin{prop}\label{cgroup}
Let $\HH = \HH(G; \g_0, \g_1, \g_2)$ be a regular hypermap with simple underlying Levi graph. Then $\lg \g_1,\g_2\rg \cap \lg\g_0,\g_2\rg=\lg\g_2\rg$.
\end{prop}

\demo
Since the Levi graph is simple, a hypervertex and a hyperedge are incident with at most 2 flags. So $|\lg \g_1,\g_2\rg \cap \lg\g_0,\g_2\rg|\le 2$. It follows that  $\lg \g_1,\g_2\rg \cap \lg\g_0,\g_2\rg=\lg\g_2\rg$.

\qed

As we know, $G$ acts faithfully on hypervertices of the regular hypermap $\HH$ if and only if the subgroup $\lg \g_1,\g_2\rg$ has a trivial core in $G$. The core of $H =\lg \g_1,\g_2\rg $ in $G$ is closely related to occurrence of multiple adjacencies in the underlying hypergraph of the hypermap $\HH$. If each $m$ adjacent vertices is
joined by exactly $t$ distinct hyperedges we will say that the hypergraph has {\it hyperedge-multiplicity} $t$; if $t = 1$ the hypergraph is {\it simple}. In what follows we will assume that our regular hypermaps are not cycles in a sphere, in particular, this means that for the valence $k$ of a hypervertex we may assume that $k\ge 3$ when $m=2$.

Let $K=\lg \g_2,\g_0\rg$ be a hyperedge stabiliser and $[HK:H]=\{Hx \di x \in K\}$. Then $K$ is a hyperedge incident with the hypervertices in $[HK : H] $. To determine the hyperedge-multiplicity of $\HH$ we may (by regularity) just count the number
of hyperedges incident with all the hypervertices in $[HK : H]$. We see that $Hx\cap Kg \ne \emptyset$ for all $x \in K$ if and only if $g\in \bigcap\limits_{x\in K}(KH^x)$. Therefore, the hyperedges incident with all the hypervertices in $[HK : H]$ are in a 1-1 correspondence with cosets of the subgroup $K$ in $ \bigcap\limits_{x\in K}(KH^x)$. Thus, for arbitrary hypermaps, we have the following proposition.

\begin{prop} \label{muti}
 The underlying hypergraph of the regular hypermap $\HH$ associated with a group $G$ with
presentation {\rm (\ref{XX})} has edge-multiplicity $|\bigcap\limits_{x\in K}(KH^x) : K|$ where $H =\lg \g_1,\g_2\rg $ and $K=\lg \g_2,\g_0\rg$.
\end{prop}

\begin{theo}
Let $\HH=\HH(G;\g_0, \g_1, \g_2)$ be a regular hypermap, $H =\lg \g_1,\g_2\rg $ be a hypervertex stabiliser and $K=\lg \g_2,\g_0\rg$ be a hyperedge stabiliser. Suppose $|\g_1\g_2|=k\ge 3$ when $|\g_2\g_0|=m=2$.
Then,
\begin{itemize}
\item[$(i)$]
 the underlying hypergraph of $\HH$ is simple if and only if $K=\bigcap\limits_{x\in K}(KH^x)${\rm ;}
 \item[$(ii)$]
  if $K=\bigcap\limits_{x\in K}(KH^x)$, then $\Core_G(H)=1$.
\end{itemize}
\end{theo}

\demo
By Proposition \ref{muti}, we obtain (i). Now we proof (ii).

 Suppose that $N=\Core_G(H)$. Assume that $N\ne 1$. Then $K=\bigcap\limits_{x\in K}(KH^x)$ implies $N\le K$. And so $N\le H\cap K$. Note that $H\cap K = \lg \g_2 \rg$ from Proposition \ref{cgroup}. Hence $N=\lg \g_2\rg$. It follows that
$$ G=\langle r_2, r_1, r_0 \di r_2^2=r_1^2=r_0^2=1, (r_0r_1)^n=(r_2r_1)^2=(r_2r_0)^2=1\rangle \cong D_{2n}\times \mathbb{Z}_2,$$
which contradicts $k\ge 3$ when $m=2$. Hence $N=1$, that is, $\Core_G(H)=1$, as desired.

\qed

\begin{rem}Let $G=\lg \g_0, \g_1, \g_2\rg, H =\lg \g_1,\g_2\rg $ and $K=\lg \g_2,\g_0\rg$ where $|\g_0|=|\g_1|=|\g_2|=2$.
If $|H\cap H^{\g_0}|=2$ and $H \cap K =\lg \g_2 \rg$, then $|\bigcap\limits_{x\in K}(KH^x) : K|$ may be great than 1; if $|H \cap K|=4>2$, then $|\bigcap\limits_{x\in K}(KH^x) : K|$ may be equal to 1. For example, $G=PSL(2,11)$.
\end{rem}

\begin{cor}
Let $\HH=\HH(G;\g_0, \g_1, \g_2)$ be a regular hypermap with a underlying simple hypergraph. Suppose $|\g_1\g_2|=k\ge 3$ when $|\g_2\g_0|=m=2$. Then the automorphism group of $\HH$ has a faithful action on its  hypervertices.
\end{cor}

Now we obtain the following classification algorithm.
\begin{alm} \label{alm}{\rm
 To classify the  regular embeddings of simple hypergraphs with given order $\omega$, we need to do the following two steps:
\begin{itemize}
\item[$(i)$]
Up to the group isomorphism, find all the transitive subgroups $G$ in the symmetric group $S_\omega$ of degree $\omega$ on a set $V$, such that
$$G = \langle \g_0, \g_1, \g_2 \di \g_0^2=\g_1^2=\g_2^2=(\g_1\g_2)^k=(\g_2\g_0)^m=(\g_0\g_1)^n= 1\rangle,$$
$$G_v=\langle \g_1,\g_2\rangle, K=\lg\g_2,\g_0 \rg,$$
$$G_v \cap K = \lg \g_2 \rg, K=\bigcap\limits_{x\in K}(KG_v^x).$$
\item[$(ii)$] For each group $G$ in $(i)$, determine the representatives of the orbits of $\Aut(G)$ on the set of triples $(\g_0,\g_1,\g_2)$, satisfying the
conditions in $(i)$. Then we obtain all the non-isomorphic  regular simple hypermaps $\HH(G; \g_0,\g_1,\g_2)$, whose numbers of hypervertices and flags are $\omega$ and $|G|$, respectively.
\end{itemize}
}
\end{alm}

As an application of Algorithm \ref{alm}, we determine the automorphism groups of regular embedding of simple hypergraphs
with prime square order.

\section{Some Results}

Throughout the paper, we fix the notation for two finite fields: $\B{F}_p^*=\lg \th \rg ,$  $\B{F}_{p^2}=\B{F}_{p}(\a), $ where $\a^2=\th$.
We shall adopt matrix notation for $\GL(n, p)$,
denote a matrix $x=(a_{ij})_{n\times n}$ by $x=||a_{11},a_{12},\cdots,a_{1n};\cdots;a_{n1},a_{n2}\cdots,a_{nn}||$.
In an affine geometry $\AG(V)$,  for any row vector $\a$ in $V,$  by $t_\a $ we denote the translation corresponding to $\a $
and  by $T$ the translation subgroup of $\AGL(n,p)$ so that $\AGL(n, p)\cong T\!\rtimes\!\GL(n, p)$.
Then for any $t_\a \in T\leq \AGL(n,p)$ and any $g\in \GL (n, p)\le \AGL(n,p),$  we have $g^{-1}t_\a g=(t_\a)^g=t_{\a g}$.
Finally, by $(r,s)$ and $[r,s],$ we denote the greatest common divisor and the least common multiple  of two positive integers $r$ and $s,$ respectively.  For a group $G$, by $G^k$ we denote the subgroup $\lg g^k\di g\in G\rg$ of G.
For  positive integers $k$, by $\pi(k)$  we denote the set of all the primes which are divisors of $k$.

\begin{prop} {\rm \cite[Theorem 3.4, 3.5]{DMB}} \label{GL(2,p)}
For an odd prime $p$, let $H$ be a maximal subgroup of  $G=\GL(2,p)$. Then   up to conjugacy, $H$ is isomorphic  to one of the following subgroups:
\begin{enumerate}
\item[\rm{(1)}] $D\rtimes \lg b\rg$, where $D$ is the subgroup of diagonal matrices and $b=||0,1;1,0||;$
\item[\rm{(2)}] $\lg a\rg\rtimes\lg b\rg$, where $\lg a\rg$ is a Singer subgroup of $G$ and $b=||1,0;0,-1||,$ defined by $a=||e,f\th;f,e||\in G$,
where $\B{F}_p^*=\lg \th\rg$, $\B{F}_{p^2}=\B{F}_p(\a)$ for $\a^2=\th$ and $\B{F}_{p^2}^*=\lg e+f\a \rg;$
\item[\rm{(3)}] $\lg a\rg\rtimes D$, where $a=||1,1;0,1||;$
\item[\rm{(4)}] $H/\lg z\rg$ is isomorphic to $A_4\rtimes \B{Z}_{\frac{p-1}2}$ for $P\equiv 5(\mod 8);$
$S_4\rtimes \B{Z}_{\frac{p-1}2}$ for $P\equiv 1,3,7(\mod 8);$ or $A_5\rtimes \B{Z}_{\frac{p-1}2}$ for $P\equiv \pm 1(\mod 10),$ where
$z=||-1,0;0,-1||$ and   $\B{Z}_{\frac{p-1}2}=Z(G)/\lg z\rg;$
\item[\rm{(5)}]   $H/\lg z\rg$  is isomorphic to $A_4\rtimes \lg s\rg$, where $\lg s^2\rg\leq Z(G)/\lg z\rg$ for $p\equiv 1(\mod 4). $
\end{enumerate}
\end{prop}

\begin{prop} {\rm  \cite[Lemma 2.6]{DK}} \label{2k}
For an odd prime $p$, let $G=\GL(2,p)$. Then  any dihedral subgroup $H\leq G$ of order $2k$ is conjugate to one of the following subgroups.
\begin{enumerate}
\item[\rm{(1)}] $\lg a\rg\rtimes\lg b\rg$, where $a=||t,0;0,t^{-1}||$ and $b=||0,1;1,0||.$  Note that $|t|=k\ge 3$ and $k\di (p-1).$
\item[\rm{(2)}] $\lg a\rg\rtimes\lg b\rg$, where $a=||e,f\th;f,e||$ and $b=||1,0;0,-1||.$ Note that $e^2-f^2\th=1,$ $|e+f\a|=k\ge 3$ and $k\di (p+1).$
\item[\rm{(3)}] $\lg a\rg\rtimes\lg b\rg$, where $a=||i,i;0,i||$ for $i=1$ or $-1$ and $b=||1,0;0,-1||.$
Note that $k=p$ if  $i=1$, and  $k=2p$ if $i=-1.$
\item[\rm{(4)}] $\lg a\rg\rtimes\lg b\rg$, where $a=||-1,0;0,-1||$ and $b=||1,0;0,-1||.$ Note that $k=2.$
\end{enumerate}
\end{prop}

\begin{prop} \label{p^k} {\rm \cite{Gur}}\,
Let $T$ be a nonabelian simple group with a subgroup $H<T$ satisfying $|T:H|=p^k$ for $p$ a prime. Then one of the following holds:
\begin{enumerate}
\item[\rm {(1)}]  $T=A_n$ and $H=A_{n-1}$ with $n=p^k;$
\item[\rm{(2)}]  $T=\PSL(n,q)$ and $H$ is the stabilizer of a projective point or a hyperplane in $\PG(n-1,q),$ and $|T:H|=(q^n-1)/(q-1)=p^k;$
\item[\rm{(3)}]  $T=\PSL(2,11)$ and $H=A_5;$
\item[\rm{(4)}]  $T=M_{11}$ and $H=M_{10};$
\item[\rm{(5)}]  $T=M_{23}$ and $H=M_{22};$
\item[\rm{(6)}] $T=\PSU(4,2)$, $\PSp(4,3)$ and $H$ is a subgroup of index $27$.
\end{enumerate}
\end{prop}

\begin{prop}{\rm \cite[17, Satz 4.5]{Hup}}\label{NC}
Let $H$ be a subgroup of a group $G$. Then $C_G(H)$ is a normal subgroup of $N_G(H)$ and the quotient $N_G(H)/C_G(H)$ is isomorphic with a subgroup
of $\Aut(H)$.
\end{prop}

\begin{prop} {\rm \cite[11.6, 11.7]{Wie}} \label{p}
Every permutation group of prime degree $p$ is  isomorphic to either  $\B{Z}_p\rtimes \B{Z}_s$ for some $s\di (p-1)$
or an insolvable group listed in Proposition~\ref{p^k} for $k=1$.
\end{prop}

From Proposition \ref{p}, one can get the following lemma.

\begin{lem} \label{psemi} Suppose that $G$ is a permutation group of degree $p$ where $p$ is a prime,
whose point stabilizers are either cyclic groups or dihedral groups. Then  $G\cong\B{Z}_p\rtimes\B{Z}_s$ where $s\di (p-1)$.
Moreover, the point stabilizer is  a cyclic group.
\end{lem}

\begin{prop} \label{mn}  Let  $G=MN$, where $M\unlhd G$ and $N\unlhd G$, $M\cap N=1$. Then $G=M\times N$.
\end{prop}

\begin{lem} \label{mg}
Suppose that   $M$ and $Q$ are normal $p$-subgroup of $G$ with $M< Q$,  $L\leq G$ where $(|L|,p)=1$.   Suppose that  $ M\rtimes L\unlhd G$.
Then for every $g\in Q\setminus M$ and $l\in L$, we have that  $g^l= gm$ for some $m\in M$.
\end{lem}
\demo For every $g\in Q\setminus M$ and $l\in L$,  we have that  $l^g=l(g^{-1})^lg\in M\rtimes L$, which forces that $(g^{-1})^lg\in  M\rtimes L$.
It follows that $g^{-1}g^l\in M\rtimes L$. Note that $Q\unlhd G$. Then  $g^{-1}g^l\in Q$. Since $(|L|,p)=1$ and $M< Q$, we get  $g^{-1}g^l\in M$.  Thus,  $g^l=gm$ for some $m\in M$.
\qed

\section{Transitive Permutation Groups}

In this section, as an application of Algorithm \ref{alm}, in order to classify the  regular embeddings of simple hypergraphs with given order prime square, we determine all the transitive permutation groups of degree prime square, whose point stabilizers are  dihedral groups.

\begin{theo}\label{main} Suppose that $G$ is a transitive permutation group of degree $p^2$ where $p$ is a prime, whose point stabilizers are  dihedral groups. Let $T\cong\B{Z}_{p}\times \B{Z}_p$. Then $G$ is isomorphic  to one of the following groups:
\begin{enumerate}
\item[\rm{(1)}]  $G_1=S_4\cong T\rtimes\lg x,y\rg$, where  $p=2$ and $n=3,$ $x=\|1,1;1,0\|, y=\|0,1;1,0\|;$

\item[\rm{(2)}]  $G_2=T\rtimes\lg x,y\rg$, where $p\ge 3$, $n\di (p+1)$ and $n\ge 3$,  $x=||e,f\th;f,e||$ with $e^2-f^2\th=1$ and $y=||1,0;0,-1||;$

\item[\rm{(3)}]  $G_3=T\rtimes\lg x,y\rg$, where $p\ge 3$, $n\di (p-1)$and $n\ge 3$, $x=||t,0;0,t^{-1}||$ and $y=||0,1;1,0||;$

\item[\rm{(4)}]  $G_4=T\rtimes\lg x,y\rg$,  where $p\ge 3$, $n=2$, $x=||-1,0;0,-1||$ and $y=||1,0;0,-1||;$

\item[\rm{(5)}]  $G_5=T\rtimes\lg x,y\rg$,  where $p\ge 3$, $n=p$, $x=||1,1;0,1||$ and $y=||-1,0;0,1||;$

\item[\rm{(6)}]  $G_6=T\rtimes\lg x_1,x_2,y\rg$, where $p\ge 3$, $n=2p,$   $x_1=||1,1;0,1||$, $x_2=||-1,0;0,-1||$  and $y=||1,0;0,-1||;$

\item[\rm{(7)}]   $G_7=T\rtimes\lg x,y\rg$, where $p\ge 3$,  $n=p,$ $x=||1,1;0,1||$ and $y=||-1,0;0,-1||.$
\end{enumerate}
\end{theo}

Throughout this section, let  $G$ be a transitive permutation group of degree $p^2$ and  $G_v\cong\B{D}_{2n}$, where $v\in\V(\G)$.
Note that $|G|=2np^2$.  Pick up a subgroup $P\in\Syl_p(G)$.

\vskip 3mm
If $p$=2, then $n=3$ and so  $G\cong S_4$. So from now on we suppose that  $p\ge 3$ and $T\cong\B{Z}_{p}\times \B{Z}_p$.

\subsection{ $ G$  is a primitive group }

\begin{lem} \label{pri}   Suppose that $G$ is a primitive group of order $p^2$. Then
$G=T\rtimes\lg x,y\rg$ is an affine group,  and $\lg x,y\rg$ is conjugate to one of the following subgroups:
\begin{enumerate}
\item[\rm{(1)}] $x=||e,f\th;f,e||$ and $y=||1,0;0,-1||,$ where $p\ge 3$, $|x|=n\ge 3$ and $n\di (p+1);$
\item[\rm{(2)}] $x=||t,0;0,t^{-1}||$ and $y=||0,1;1,0||,$ where $p\ge 3$, $|x|=n\ge 3$ and $n\di (p-1)$.
\end{enumerate}
\end{lem}
\demo Suppose that $G$ is a primitive group of order $p^2$.
Then  by checking  O'Nan-Scott Theorem \cite[Theorem 4.1A]{DM}, we know that every  primitive group of degree $p^2$  is  an almost simple group,
or an affine group. If $G$ is an almost simple group, then  by Proposition~\ref{p^k}
we know that for an  almost simple group of degree $p^2$, its point stabilizer  cannot be a dihedral group, a contradiction.
Thus, $G$ is  an affine group, and then  $G=T\!\rtimes\!G_v\cong \B{Z}_p^2\!\rtimes\!\B{D}_{2n}$  where $G_v$ is an irreducible subgroup of $\GL(2,p)$.  From Proposition \ref{2k}, one can get $\lg x,y\rg$ is conjugate to one of the subgroups in (1) and (2). This completes the proof.
\qed

\subsection{ $ G$  is an imprimitive group }

Suppose that $G$ is an imprimitively of order $p^2$.  Let ${\cal B}$  be a complete block system of $G$.
Let $K$  be the kernel of $G$ on the block system.
If $K=1$, then $G\cong G/K\lesssim S_p$, contradicting  to $p^2\di |G|$. Thus,  $G/K  \lesssim S_p$ and  $K\neq 1$.

Considering  $K$ acts on the block $B$, where $B\in\cal{B}$. Then $K$ is transitive  on each block.
For any block $B\in{\cal B}$ and any vertex $v\in B$, we have $$(G/K)_B=G_vK/K\cong G_v/(G_v\cap K)=G_v/K_v. $$
Since $G_v\cong\B{D}_{2n}$,  $G_v/K_v$ is cyclic or dihedral. Note that $G/K\lesssim S_p$. Then $(G/K)_B$ is cyclic, that is $G_v/K_v$ is cyclic.
Since $K_v\unlhd G_v$,  we have that
\begin{enumerate}
\item[\rm {(1)}]  if $K_v$  is a cyclic group, then $K_v\cong\B{Z}_n$ and $G/K\cong\B{Z}_p\rtimes\B{Z}_2;$
\item[\rm{(2)}] if $K_v$ is a dihedral group, then either $K_v\cong\B{D}_{2n}$ and $G/K\cong\B{Z}_p;$ or $K_v\cong\B{D}_{n}$ and
$G/K\cong\B{Z}_p\rtimes\B{Z}_2$.
\end{enumerate}

Let $K_0$ be the  kernel of $K$ on the block $B$.
Then $K/K_0\lesssim S_p$. Set $G_v=\lg x,y\rg\cong\B{D}_{2n}$, where  $|x|=n$. Take $P\in\Syl(G)$.

\begin{lem} \label{k0=1}   Suppose that $K_0=1$. Then $n=2$ and $G=T\rtimes\lg x,y\rg$ is an affine group,
where $x=||-1,0;0,-1||$ and $y=||1,0;0,-1||$.
\end{lem}
\demo Suppose that $K_0=1$.  Then $K\lesssim S_p$. Note that $G/K\lesssim S_p$.
Then by Lemma \ref{psemi}, we have $K\cong\B{Z}_p\rtimes\B{Z}_n$ where $n\di (p-1)$,  and $G/K\cong\B{Z}_p\rtimes\B{Z}_2$.
Let $P\in\Syl(G)$. Then   $PK/K\unlhd G/K$ and so $PK\unlhd G$. Since $|PK|=\frac{|P||K|}{|P\cap K|}=p^2n$, we get $P\unlhd PK$.
Thus $P\unlhd G$ and $|P|=p^2$. Moreover, $G=P\rtimes G_v$.

Suppose that $P\cong\B{Z}_{p^2}$. Then we set $P=\lg a\rg$  and   $K=\lg a^p\rg\rtimes\lg x\rg$, where $|x|=n$.
By Lemma \ref{mg}, we set $a^x=a^{1+kp}$. Thus,  $(a^p)^x=a^p$, which implies that $K=\lg a^p\rg\times\lg x\rg$.
Since $(|x|,p)=1$, we get that $\lg x\rg$ is a   characteristic subgroup of $K$.
Therefore, $\lg x\rg\unlhd G$,  which implies that it fixes all vertices, a contradiction.

Suppose that $P\cong\B{Z}_p^2$. Then $P\leq C_G(P)$. If $P< C_G(P)$,  then there exists an element  $g$ in  $G_v$ such that $g\in  C_G(P)$.
Noting that  $P$ is a transitive group, we know that  $g$ fixes all vertices, a contradiction. Thus $P=C_G(P)$.
By Proposition \ref{NC}, we have that $$\B{D}_{2n}\cong G/P=N_G(P)/C_G(P)\leq  \Aut(P)\cong \GL(2,p).$$
Note that $G$ is an imprimitively group and $|G|=2np^2$ where $n \di (p-1)$. Then by  checking the Proposition \ref{2k}
we know that $n=2$ and $G=T\rtimes\lg x,y\rg$  is an affine group, where $x=||-1,0;0,-1||$ and $y=||1,0;0,-1||$.

\begin{lem} \label{kuc}   Suppose that $K_0\neq 1$ and $K_v\cong\B{Z}_n$. Then $n=p$ and
$$G=\lg a,b,x\di a^p=b^p=c^p=x^2=1, [a,b]=c, [a,c]=[b,c]=1, a^x=a^{-1}, b^{x}=b^{-1}\rg. $$
\end{lem}
\demo Suppose that $K_0\neq 1$ and $K_v\cong\B{Z}_n$.  Then $G/K\cong\B{Z}_p\rtimes\B{Z}_2$,
$K_0\leq K_v$ and by Lemma \ref{psemi}, $K/K_0\cong\B{Z}_p\rtimes\B{Z}_s$, where $s\di (p-1)$ and $s\di n$.
Take $M\in\Syl_p(K)$ and $P\in\Syl_p(G)$. In what follows, we divide our proof into four steps.

\vskip 3mm
{\it Step 1:}   {\it Show   $M\unlhd K$ .}
\vskip 3mm
Note that $n\leq p^2-1$. Then  $|M|=p$ or $p^2$,  $MK_0/K_0\unlhd K/K_0$ and so $MK_0\unlhd K$.
Since $K_v\cong\B{Z}_n$ and $K_0\leq K_v$, we get $K_v\leq C_K(K_0)\unlhd N_K(K_0)=K$.
Since $|K|=pn$, we  have that $C_K(K_0)=K_v$ or $C_K(K_0)=K$.

Suppose that $C_K(K_0)=K_v$. Then $K_v\unlhd N_K(K_0)=K$.
If $|\pi(n)\setminus \{p\}|>1$,  we take $l\in\pi(n)\setminus \{p\}$ and $L\in\Syl_l(K_v)$.
Then  $L\in\Syl_l(K)$ and  $L$ is a cyclic characteristic subgroup of $K_v$, which forces that  $L\unlhd K$.
This implies that   $L$  is a characteristic subgroup of $K$. It follows that $L\unlhd G$,  which implies that  it fixes all vertices, a contradiction. Therefore, $|\pi(n)\setminus \{p\}|=0$, that is $K_v\cong\B{Z}_p$, and then  $|K|=p^2$, contradicting to $C_K(K_0)=K_v$.

Suppose that   $C_K(K_0)=K$. Then  $K_0\leq Z(K)$, which implies that $MK_0$ is an  abelian group.
Since $M\in\Syl_p(MK_0)$, we get that  $M$ is a characteristic subgroup of $MK_0$, which forces that $M\unlhd K$.
Therefore, we set  $K=M\rtimes\lg y\rg$, where $y\in K_v$ and $(|y|,p)=1$.

\vskip 3mm
{\it Step 2:}  {\it Show $ K=M\cong\B{Z}_p^2$.}
\vskip 3mm
Now, we know that $M\cong\B{Z}_p;$ or $M\cong\B{Z}_{p^2};$ or $M\cong\B{Z}_p^2$.

Suppose that   $M\cong\B{Z}_p$. Then $M\nleq Z(K)$. Otherwise, $K$ is a cyclic group, a contradiction.  Thus  $(|Z(K)|,p)=1$.
Note that $(|K_0|,p)=1$ and $K_v\cong\B{Z}_n$.  Then by Proposition \ref{mn} we know that  $K_0\leq Z(K)$.
Since $K=M\rtimes\lg y\rg$,  we get $Z(K)\leq \lg y\rg$. This implies that  $K_0$ is a  characteristic subgroup of $K$.
Since $K\unlhd G$, we get  $K_0\unlhd G$,  which implies that it fixes all vertices. It follows that $K_0=1$,   a contradiction.

Suppose that  $M\cong\B{Z}_{p^2}$. Then $M^p$ is a characteristic subgroup of $K$. Thus,  $M^p\unlhd G$,
which implies that it fixes all vertices. It follows that $M^p=1$,  a contradiction.

Suppose that   $M\cong\B{Z}_p^2$. Then we  set $M=\lg a,b\rg$ and  $\lg b\rg\leq  K_0$. This implies  that $\lg b\rg\leq Z(K)$.
If  $|y|\ge 2$, then $a\not\in Z(K)$. Otherwise, $M\leq Z(K)$ and so $K$ is an abelian group, a contradiction.
Thus  $ Z(K)\cap M=\lg b\rg$,  which implies that  $\lg b\rg$ is  a characteristic subgroup of $Z(K)$.
It follows that  $\lg b\rg$ is  a characteristic subgroup of $K$.
Since $K\unlhd G$, we get  $\lg b\rg\unlhd G$, which implies that it fixes all vertices,  and so $\lg b \rg=1$, contradicting to $|b|=p$.
Therefore, $y=1$ and $ K=M\cong\B{Z}_p^2$.   Moreover, $n=p$.
Note that $G/K\cong\B{Z}_p\rtimes\B{Z}_2$ and $ K\cong\B{Z}_p^2$.  Then $P\unlhd G$ and we  set $G=P\rtimes \B{Z}_2$.

\vskip 3mm
{\it Step 3:}   {\it Show   $P\cong(\B{Z}_{p}\times \B{Z}_p)\rtimes\B{Z}_p$ .}
\vskip 3mm

Note that  $P$ is a transitive group. Then either $P\cong\B{Z}_{p^2}\rtimes\B{Z}_p;$ or $P\cong(\B{Z}_{p}\times \B{Z}_p)\rtimes\B{Z}_p$.

Suppose that $P\cong\B{Z}_{p^2}\rtimes\B{Z}_p$. Then we set $$G=\lg a,b,x\di a^{p^2}=b^p=x^2=1, a^b=a^{1+p}, a^x=a^ib^j,b^x=a^{pk}b^{t}\rg,$$
where  $i\in\B{Z}_{p^2}^*$, $t\in\B{Z}_p^*$ and $x\in G_v$.
Set $P_v=\lg a^{pm}b\rg$ for some $m$,  such that $(a^{pm}b)^x=a^{-pm}b^{-1}$.
Then $$a^{-pm}b^{-1}=(a^{pm}b)^x=a^{pmi+pk}b^t.$$  This implies that $t\equiv -1(\mod p)$.
By $a^b=a^{1+p}$, we get  $$a^{i+pi}b^j=(a^{1+p})^x=(a^b)^x=(a^ib^j)^{a^{pk}b^{-1}}=a^{i-pi}b^j,$$
which  implies that $pi\equiv -pi(\mod p^2)$. It follows that $i\equiv 0(\mod p)$, a contradiction.
Therefore, $P\cong(\B{Z}_{p}\times \B{Z}_p)\rtimes\B{Z}_p$.

\vskip 3mm
{\it Step 4:}  {\it Determination of the groups $ G$. }
\vskip 3mm

Suppose that  $$G=\lg a,b,x\di a^p=b^p=c^p=x^2=1, [a,b]=c,[a,c]=[b,c]=1, a^x=a^ib^jc^k, b^x=b^{-1}\rg, $$ where $i\in\B{Z}_p^*$ and $P_v=\lg b\rg$.
Then  $c^x=[a^ib^jc^k,b^{-1}]=c^{-i}$,  $K=\lg b,c\rg$ and $G/K=\lg \oa,\ox\rg$.
Since $G/K\cong\B{Z}_p\rtimes\B{Z}_2$, we get $i\equiv -1(\mod p)$ and so $a^x=a^{-1}b^jc^k$.
By $x^2=1$, we get that  $$a=a^{x^2}=(a^{-1}b^jc^k)^x=(a^{-1}b^jc^k)^{-1}b^{-j}c^{k}=b^{-j}ab^{-j}=ab^{-2j}c^j,$$
which implies that $-2j\equiv 0(\mod p)$. It follows that $j\equiv 0(\mod p)$. Thus, $a^x=a^{-1}c^k$ and $c^x=c$.
Let $a'=ac^{-\frac{k}2}.$ Then $[a',b]=[ac^{-\frac{k}2},b]=[a,b]=c$ and
$a'^x=(ac^{-\frac{k}2})^x=a^{-1}c^kc^{-\frac{k}2}=a^{-1}c^{\frac {k}2}=a'^{-1}$. Therefore,  we set  $a^x=a^{-1}$ and $b^x=b^{-1}$.
\qed

\begin{lem} \label{kud1}  Suppose that $K_0\neq 1$ and $K_v\cong\B{D}_{2m}$, where $m=n$ or $2m=n$.
Then  $K\cong\B{Z}_p^2\rtimes\B{Z}_2$ or $K\cong\B{Z}_p^2\rtimes\B{D}_4$.
\end{lem}
\demo
Suppose that $K_0\neq 1$ and $K_v\cong\B{D}_{2m}$.  Then $G/K\cong\B{Z}_p$ for $m=n$; or $G/K\cong\B{Z}_p\rtimes\B{Z}_2$ for $2m=n$,
$K_0\unlhd  K_v$ and $K/K_0\cong\B{Z}_p\rtimes\B{Z}_s$, where $s=1$ or $s=2$. Take $M\in\Syl_p(K)$.
In what follows, we divide our proof into three steps.

\vskip 3mm
{\it Step 1:}  {\it Show  $\pi(m)\subseteq \{p,2\}$.}
\vskip 3mm

Suppose that  $|\pi(m)\setminus \{p,2\}|>1$.
Take $q\in \pi(m)\setminus \{p,2\}$ and let $Q\in\Syl_q(K_v)$.
Since $|G/K|\di 2p$ and $|K/K_0|\di 2$, we get $Q\in\Syl_q(K_0)$  and  $Q\in\Syl_q(K)$, which implies that
$Q$ is a cyclic characteristic subgroup of $K_0$. Since $K_0\unlhd K$,  we get  $Q\unlhd K$, which forces that
$Q$ is a characteristic subgroup of $K$. It follows that  $Q\unlhd G$.
Noting that  $Q\leq G_v$,   we get  $Q$ fixes all vertices and so $Q=1$,  a contradiction.
Therefore, $ |\pi(m)\setminus \{p,2\}|=0$, that is $\pi(m)\subseteq \{p,2\}$.  Moreover, $\pi(|K_0|)\subseteq \{p,2\}$.

\vskip 3mm
{\it Step 2:}  {\it Show  $M\unlhd K$ and $M\cong\B{Z}_p^2$.}
\vskip 3mm

 Note that $n\leq p^2-1$. Then  $|M|=p$ or $p^2$.

Suppose that $p\not\in\pi (m)$.  Then  $\pi(m)=2$, $|M|=p$ and $K_v\in\Syl_2(K)$.  This implies that $K_0$ is a normal $2$-subgroup of $K$.
Thus, $K_0\leq {\rm O}_2(K)$. Note that $|K_v/K_0|\di 2 $. Then ${\rm O}_2(K)=K_0$ or ${\rm O}_2(K)=K_v$.
If  ${\rm O}_2(K)=K_v$, then $K_v$ is a   characteristic subgroup of $K$,  and thus  $K_v\unlhd G$, a contradiction;
if ${\rm O}_2(K)=K_0$, then $K_0$ is a   characteristic subgroup of $K$, and thus $K_0\unlhd G$, a contradiction.
Therefore, $p\in \pi(m)$ and $|M|=p^2$.

Let $Q\in \Syl_p(K_v)$. Then $Q\in\Syl_p(K_0)$ and $Q$ is a  characteristic subgroup of $K_0$.
This implies that $Q\unlhd K$. It follows that $Q\leq {\rm O}_p(K)$.
Since $|M|=p^2$, we get ${\rm O}_p(K)=Q$ or $M$.
If  ${\rm O}_p(K)=Q$, then $Q$ is a   characteristic subgroup of $K$,  and thus  $Q\unlhd G$, a contradiction.
Thus,  ${\rm O}_p(K)=M$  is a   characteristic subgroup of $K$,  which  implies that $M\unlhd K$.

Noting that $|M|=p^2$, we have either  $M\cong\B{Z}_p^2;$ or $M\cong\B{Z}_{p^2}$.
If   $M\cong \B{Z}_{p^2}$, then  $M^p$ is a   characteristic subgroup of   $K$, which forces that  $M^p\unlhd G$, a contradiction.
Hence, $M\cong\B{Z}_p^2$.

\vskip 3mm
{\it Step 3:}  {\it Either $K\cong\B{Z}_p^2\rtimes\B{Z}_2,$  or $K\cong\B{Z}_p^2\rtimes\B{D}_4$. }
\vskip 3mm
First, suppose that   $2\not\in\pi (m)$.  Then $m=p$  and so  $|K|=2p^2$.  Note that $K_v\cong\B{D}_{2p}$.  Then $K\cong\B{Z}_p^2\rtimes\B{Z}_2$.

Next, suppose that  $2\in \pi(m)$. Then we set $m=2^tp$ where $t\ge 1$, and assume that
$$K_v=\lg b,x,y\di b^p=x^{2^t}=y^2=1,b^x=b,b^y=b^{-1},x^y=x^{-1} \rg\,\,{\rm and}\,\,  M=\lg a,b\rg.$$
Then $K=\lg a,b,x,y \rg$.  Since $K/K_0\cong\B{Z}_p\rtimes\B{Z}_s$ where $s=1$ or $2$, we have  either $\lg b,x\rg \leq K_0$,
or $\lg b,x^2\rg\leq K_0$.

Suppose that $\lg b,x\rg\leq K_0$. Then $\lg b,x\rg$ is a   characteristic subgroup of $K_0$, which implies that $\lg x\rg$ is a   characteristic subgroup of $K_0$. It follows that $\lg x\rg\unlhd K$.  By Proposition \ref{mn}, we get that $\lg a,b,x\rg$ is an abelian characteristic subgroup of $K$,
which forces that $\lg x\rg$ is   a   characteristic subgroup of $K$.  It follows that $\lg x\rg\unlhd G$, a contradiction.

Suppose that $\lg b,x^2\rg\leq K_0$. Then $\lg b,x^2\rg$ is a   characteristic subgroup of $K_0$, which implies that $\lg x^2\rg$ is a    characteristic subgroup of $K_0$. It follows that $\lg x^2\rg\unlhd K$.  By Proposition \ref{mn}, we get that $\lg a,b,x^2\rg$ is an abelian characteristic subgroup of $K$, which forces that $\lg x^2\rg$ is  a   characteristic subgroup of $K$. It follows that $\lg x^2\rg\unlhd G$.
Therefore, $x^2=1$, $m=2p$, $K_0\cong\B{Z}_p\rtimes\B{Z}_2$,  $K\cong\B{Z}_p^2\rtimes\B{D}_4$ and $K/K_0\cong\B{Z}_p\rtimes\B{Z}_2$.
 \qed

\vskip 3mm
By  Lemma \ref{kud1}, we know that $K\cong\B{Z}_p^2\rtimes\B{Z}_2$ or $K\cong\B{Z}_p^2\rtimes\B{D}_4$.   Next, we show that  $K\not\cong\B{Z}_p^2\rtimes\B{D}_4$.

\begin{lem} \label{kud2}   $K\not\cong\B{Z}_p^2\rtimes\B{D}_4$. Moreover, $K\cong\B{Z}_p^2\rtimes\B{Z}_2$ and $K_v\cong\B{D}_{2p}$.
\end{lem}
\demo
Suppose that  $K\cong\B{Z}_p^2\rtimes\B{D}_4$.    Then $m=2p$, $K_0\cong\B{Z}_p\rtimes\B{Z}_2$ and  $K/K_0\cong\B{Z}_p\rtimes\B{Z}_2$.
Note that $G/K\cong\B{Z}_p\rtimes\B{Z}_s$, where $s=1$ or $2$.  Then $|G|=8p^3$ or $4p^3$.
Take $P\in\Syl_p(G)$. Then $|P|=p^3$.  In what follows, we divide our proof into three steps.

\vskip 3mm
{\it Step 1:}  {\it $C_G(M)=M$.}
\vskip 3mm
By $$ C_K(M)/M\unlhd N_K(M)/M=K/M\cong\B{D}_4\,\,{\rm and}\,\, K\cong\B{Z}_p^2\rtimes\B{D}_4,$$
we have either  $C_K(M)=M$ or $C_K(M)=M\times \B{Z}_2$.
If  $C_K(M)=M\times \B{Z}_2$, then $\B{Z}_2\leq Z(K)$, which implies  that it is a  characteristic subgroup of $K$,
it follows that  it is a  normal subgroup of $G$,  a contradiction.
Thus,   $C_K(M)=M$.

Note that $M\unlhd G$. Then $C_G(M)K\unlhd G$. Since $M=C_K(M)\leq C_G(M)$, we get that  $M\leq C_G(M)$  and
$$C_G(M)/M=C_G(M)/C_K(M)\cong C_G(M)K/K\unlhd G/K\cong\B{Z}_p\rtimes\B{Z}_s,$$ where $s=1$ or $2$.
Then  $P\leq C_G(M)$ or $C_G(M)=M$. If $P\leq C_G(M)$, then $P$ is an abelian group,  contradicting to $P$ is a transitive group.
Hence,  $C_G(M)=M$.

\vskip 3mm
{\it Step 2:}  {\it $P\unlhd G$.}
\vskip 3mm
By Proposition \ref{NC}, we get $$\B{Z}_p\cong P/M\leq  G/M=N_G(M)/C_G(M)\lesssim\Aut(M)\cong\GL(2,p).$$
By Proposition  \ref{GL(2,p)},
we have  either $G/M \lesssim \lg ||1,1;0,1||\rg\rtimes D$,  where $D$ is  the subgroup of diagonal matrices;  or $(G/M)/\lg z\rg \lesssim S_4$ where $z=||-1,0;0,-1||$ and $p=3$. Suppose that   $(G/M)/\lg z\rg \lesssim S_4$, where $p=3$.  Then $(G/M)/\lg z\rg\cong\B{D}_6$ or $A_4$.
Then  $\B{Z}_2\cong (K/M)/\lg z\rg\unlhd (G/M)/\lg z\rg$, that is $(G/M)/\lg z\rg$ has normal  subgroup of order $2$,  a contradiction.
Therefore, $G/M \lesssim \lg ||1,1;0,1||\rg\rtimes D$ and thus $P\unlhd G$.

\vskip 3mm
{\it Step 3:}  {\it $G$  does not exist.}
\vskip 3mm
Note that $\B{D}_4\cong K/M\unlhd G/M$ and $\B{Z}_p\cong P/M\unlhd G/M$.
Then by Proposition \ref{mn}, we get that $ P/M\times K/M\cong\B{Z}_p\times\B{D}_4$ is an abelian normal subgroup of $G/M$.
By  $G/M\lesssim \lg ||1,1;0,1||\rg\rtimes D$,  we know that   $\lg ||1,1;0,1||\rg\rtimes D$ has no such  abelian subgroup, a contradiction.
\qed

\begin{lem} \label{kud}   Suppose that $K_0\neq 1$ and $K_v\cong\B{D}_{2p}$, where $n=p$ or $2p$.
Then  $G$ is isomorphism  to one of the following groups:
\begin{enumerate}
\item[\rm {(1)}]  $T\rtimes\lg x,y\rg$, where  $n=p,$ $x=||1,1;0,1||$ and $y=||-1,0;0,-1||;$

\item[\rm{(2)}] $T\rtimes\lg x_1,x_2,y\rg$, where  $n=2p,$   $x_1=||1,1;0,1||$, $x_2=||-1,0;0,-1||$  and $y=||1,0;0,-1||$.
\end{enumerate}
\end{lem}
\demo Now, we known that $K\cong \B{Z}_p^2\rtimes \B{Z}_2$ and $G/K\cong\B{Z}_p\rtimes\B{Z}_s$,  where $s=1$ or $2$.
Take $M\in\Syl_p(K)$. Then  $C_K(M)=M$ and $M\unlhd G$.
Take $P\in\Syl_p(G)$. Then $PK/K\unlhd G/K$ and so $PK\unlhd G$. Since $|PK|=\frac{|P||K|}{|P\cap K|}=2p^3$, we get $P\unlhd PK$,
which forces that $P\unlhd G$. Thus, we have either $G=P\rtimes\B{Z}_2;$ or $G=P\rtimes \B{D}_4$.   Moreover, $P\cong\B{Z}_{p^2}\rtimes\B{Z}_p$,
or $P\cong(\B{Z}_{p}\times \B{Z}_p)\rtimes\B{Z}_p$.

\vskip 3mm
{\it Step 1:}  {\it Show $  P\cong(\B{Z}_{p}\times \B{Z}_p)\rtimes\B{Z}_p.$ }
\vskip 3mm
Suppose that $P\cong\B{Z}_{p^2}\rtimes\B{Z}_p$. Then, we set $$P=\lg a,b\di a^{p^2}=b^p=1, a^b=a^{1+p}\rg\,\,{\rm and }\,\,P_v=\lg a^{pl}b\rg,$$
where $t\in\B{Z}_p$.
Since $G_v\cong \B{D}_{2p}$ or $\B{D}_{4p}$,  exists an involution $y\in G_v$ such that $(a^{pl}b)^y=(a^{pl}b)^{-1}$.
Now, assume that $a^y=a^ib^j$ and $b^y=a^{pt}b^k$, where $i\in\B{Z}_p^*$ and $k\in\B{Z}_p^*$.
Then $$a^{-pl}b^{-1}=(a^{pl}b)^y=a^{pli+pt}b^k.$$  This implies that $k\equiv -1(\mod p)$.
By $a^b=a^{1+p}$, we get  $$a^{i+pi}b^j=(a^{1+p})^y=(a^b)^y=(a^ib^j)^{a^{pt}b^{-1}}=a^{i-pi}b^j,$$
which  implies that $pi\equiv -pi(\mod p^2)$. It follows that $i\equiv 0(\mod p)$, a contradiction.
Therefore, $P\cong(\B{Z}_{p}\times \B{Z}_p)\rtimes\B{Z}_p$.

\vskip 3mm
{\it Step 2:}  {\it Determination of the groups $ G$. }
\vskip 3mm

By $$C_G(M)/M=C_G(M)/C_K(M)\cong C_G(M)K/K\unlhd G/K\cong\B{Z}_p\rtimes\B{Z}_s,$$ we get that $P\leq C_G(M)$ or $C_G(M)=M$.
If  $P\leq C_G(M)$, then $P$ is an abelian group, contradicting to $P$ is a transitive group. Thus,  $C_G(M)=M$.
By Proposition \ref{NC}, we have that $$P/M\unlhd G/M=N_G(M)/C_G(M)\lesssim\Aut(M)\cong\GL(2,p).$$
By Proposition  \ref{GL(2,p)}, we get that  $G/M \lesssim \lg ||1,1;0,1||\rg\rtimes D$, where $D$ is  the subgroup of diagonal matrices.

Suppose that $G=P\rtimes\B{Z}_2$. Then $G/K\cong\B{Z}_p$.  Note that  $\B{Z}_2\cong K/M\unlhd G/M$. Then   by  Proposition \ref{mn},
we get that $G/M=P/M\times K/M\cong\B{Z}_p\times\B{Z}_2$.  By Proposition \ref{GL(2,p)}, we have that $G/M=\lg x,y\rg$, where $x=||1,1;0,1||$ and $y=||-1,0;0,-1||$. Therefore,  $G=T\rtimes\lg x,y\rg$.

Suppose that $G=P\rtimes\B{D}_4$. Then $G/K\cong\B{Z}_p\times\B{Z}_2$. Note  that  $\B{Z}_2\cong K/M\unlhd G/M$.
Then   by  Proposition \ref{mn},  we get that $PK/M=P/M\times K/M\cong\B{Z}_p\times\B{Z}_2$.
By Proposition \ref{GL(2,p)}, we have that  $G/M=\lg x_1,x_2,y\rg$, where $x_1=||1,1;0,1||$, $x_2=||-1,0;0,-1||$ and $y=||1,0;0,-1||$.
Therefore, $G=T\rtimes\lg x_1,x_2,y\rg$.
\qed

\vskip 0.3cm
\noindent {\bf Proof of Theorem \ref{main}} One can get this theorem from the above discussion.

\section{Regular Hypermaps of Order Prime Square}

In this section,  for $i\in \{1,2,\ldots, 7\}, $ $G_i$ is the group in Theorem \ref{main}. We determine the automorphism groups of
regular embedding of simple hypergraphs with prime square order.

\begin{lem}\label{nong6}
Let~$H=\langle r_1, r_2\rangle\cong D_{2p}$ with $|r_1|=|r_2|=2$ be a subgroup of~$G_7$. Suppose that  $\Core_G(H) = 1$ and $r_0\in G$, where~$|r_0|=2$. Then~$\langle r_0, r_1, r_2\rangle\ne G_7$.
\end{lem}

\demo
It is easy to see that $G_7$ has two conjugacy classes of subgroups~$D_{2p}$, and their representative elements are ~$\langle t_{(0,1)}, y\rangle$ and ~$\langle t_{(1,0)}, y\rangle$. If $H=\langle t_{(0,1)}, y\rangle$, then $\Core_G(H) \ne 1$ from $\langle t_{(0,1)} \rangle\unlhd G$, contradicting to~$\Core_G(H)=1$. So we may assume $H=\langle t_{(1,0)},y \rangle$. Since~$H$ has only one conjugacy class of involutions, we may choose~$r_2=y$. Set~$r_0=t_{(a, b)}x^iy$.  From $|r_0|=2$, we have $r_0^2=x^{2i}t_{(0,ai)}=1$, and so  $x^{2i}=1$, which implies that $i\equiv0(\mod p)$ and $r_0=t_{(a,b)}y$. In this case, $\langle r_2, r_1, r_0\rangle\leq \langle t_{(a,b)},y\rangle$ which is a proper subgroup of $G$. Therefore, $\langle r_0, r_1, r_2\rangle\ne G_7$.
\qed

\begin{theo}\label{main1} Let $G$ be the automorphism group of a regular simple hypermap of order  $p^2$ where $p$ is a prime. Let $T\cong\B{Z}_{p}\times \B{Z}_p$. Then $G$ is isomorphic  to one of the following groups:
\begin{enumerate}[{\rm(1)}]

\item  $G_1=T\rtimes\langle x, y\rangle\cong S_4$, where $p=2, x=\|1,1;1,0\|, y=\|0,1;1,0\|;$

\item  $G_2=T\rtimes\lg x,y\rg$, where $p\ge 3$, $n\di (p+1)$ and $n\ge 3$,  $x=||e,f\th;f,e||$ with $e^2-f^2\th=1$ and $y=||1,0;0,-1||;$

\item  $G_3=T\rtimes\lg x,y\rg$, where $p\ge 3$, $n\di (p-1)$and $n\ge 3$, $x=||t,0;0,t^{-1}||$ and $y=||0,1;1,0||;$

\item  $G_4=T\rtimes\lg x,y\rg$,  where $p\ge 3$, $n=2$, $x=||-1,0;0,-1||$ and $y=||1,0;0,-1||;$

\item  $G_5=T\rtimes\lg x,y\rg$,  where $p\ge 3$, $n=p$, $x=||1,1;0,1||$ and $y=||-1,0;0,1||;$

\item  $G_6=T\rtimes\lg x_1,x_2,y\rg$, where $p\ge 3$, $n=2p,$   $x_1=||1,1;0,1||$, $x_2=||-1,0;0,-1||$  and $y=||1,0;0,-1||$.
\end{enumerate}
\end{theo}

\demo
By Algorithm \ref{alm}, one can check that the following regular hypermaps are simple.
\begin{enumerate}[{\rm(1)}]
\item  $\HH_1=\HH(G_1;t_{(1,1)}y, xy, y)$;

\item  $\HH_{21}=\HH(G_2;t_{(1,0)}x^{\frac{n}{2}}y, xy, y)$ where $n$ is even;

\item  $\HH_{22}=\HH(G_2;t_{(0,1)}y, xy, y)$ where $n$ is odd;

\item  $\HH_{31}=\HH(G_3;t_{(1,1)}x^{\frac{n}{2}}y, xy, y)$ where $n$ is even;

\item  $\HH_{32}=\HH(G_3;t_{(1,-1)}y, xy, y)$ where $n$ is odd;

\item  $\HH_{4}=\HH(G_4;t_{(1,1)}x, y, x)$;

\item  $\HH_{5}=\HH(G_5;t_{(1,0)}y, xy, y)$;

\item  $\HH_{6}=\HH(G_6;t_{(1,0)}x_2y, x_1x_2y, y)$.
\end{enumerate}
Lemma \ref{nong6} implies that $G$ is not isomorphic to $G_7$. Therefore the theorem is proved.

\qed

\section*{Acknowledgements}
This work  was
supported by  the National Natural Science Foundation of China (No. 12101535)

\end{document}